%% file: rinv-qs-D4.tex
\documentclass[11pt]{amsart}

\usepackage{amscd}
\usepackage{amssymb}

\input skip-thm

\input qformsl-def

\newcommand{\Gb}{\bar{G}}

\begin{document}

\title[Rost invariant for quasi-split trialitarian groups]{The Rost invariant has zero kernel for quasi-split trialitarian groups}

\author{Skip Garibaldi}


\begin{abstract}
\emph{The Book of Involutions} includes the \emph{non}trivial parts of a proof that the kernel of the Rost invariant is zero for quasi-split trialitarian groups.  We record the missing (mechanical) details here for the convenience of readers who like the style of \emph{The Book}.

A proof in Harder's style can be found in a paper by Chernousov.
\end{abstract}

\maketitle

Recall that for a simply connected and absolutely almost simple algebraic group $G$ over a field $F$, \emph{the Rost invariant} is a function $H^1(F, G) \ra H^3(F, \QZt)$.\footnote{Strictly speaking, we abuse language and write ``the Rost invariant" when we really mean only the portion of it over a given field $F$. An example of the distinction is that the Rost invariant in our sense is totally uninteresting in the case $F = \C$ whereas the actual Rost invariant from \cite{GMS} or \cite[\S31]{KMRT} remains interesting.}  If $G$ is of type $\iiiD$ or $\viD$---i.e., \emph{is trialitarian}---and $F$ has characteristic $\ne 2, 3$ then the image of the Rost invariant lies in the Galois cohomology group $H^3(F, \Zm6)$; see \cite[p.~149]{GMS}.

The purpose of this note is to prove the following well-known result.

\begin{prop*}
Let $G$ be a simply connected trialitarian group over a field $F$.  If $G$ is quasi-split, then the kernel of the Rost invariant $H^1(F, G) \ra H^3(F, \Q/\Z(2))$ is zero.
\end{prop*}

You can find the important part of a proof in 40.16 of \emph{The Book of Involutions}, \cite{KMRT}.  A few technical steps are required to connect the results in  \emph{The Book} with the proposition.  They are standard for the expert familiar with the details but may be inobvious for the casual consumer, so we record them here for convenience of reference.  There is nothing interesting here.

Alternatively, \cite[Th.~6.14]{Ch:rinv} gives a proof in the style of Harder, referring to \cite{PlatRap} for details.

\section*{Preparation}

Here is a baby analogue of the proposition.  We write $F_4$ for a split simple algebraic group of that type.  It is necessarily simply connected.

\begin{lem*}
The kernel of the Rost invariant $H^1(F, F_4) \ra H^3(F, \QZt)$ is zero.
\end{lem*}

\begin{proof}
We first assume that $F$ has characteristic $\ne 2, 3$ so that we may apply results from \cite{KMRT}.
The set $H^1(F, F_4)$ classifies Albert $F$-algebras.  If such an algebra $J$ is in the kernel of the Rost invariant, then it is a first Tits construction because the mod 2 part is zero \cite[40.5]{KMRT} and it has zero divisors because the mod 3 part is zero \cite[40.8(2)]{KMRT}.  So $J$ is split by \cite[p.~416, Th.~20]{Jac:J}.

The case where $F$ has characteristic 2 or 3 follows from the characteristic zero case by the main result of \cite{Gille:inv}.
\end{proof}

\section*{The proof}
\begin{proof}[Proof of the proposition]
We assume that $F$ has characteristic zero; as in the proof of the lemma, this suffices.

The group $G$ is the identity component of the automorphism group of a twisted Hurwitz composition $\G(C, L)$ where $C$ is a split octonion algebra and $L$ is a cubic \'etale $F$-algebra.  The Springer Decomposition as in \cite[\S38.A]{KMRT} gives an injection
\[
G \injects \aut(\G(C, L)) \injects \aut(J(\G(C, L))),
\]
where the last group is of type $F_4$ \cite[38.7]{KMRT}.  In fact, $\aut(J(\G(C, L)))$ is split: it contains the rank 2 split subgroup $\aut(C)$ of type $G_2$, and a group of type $F_4$ with rank $\ge 2$ is split by Tits's classification \cite[p.~60]{Ti:Cl}.

Fix $\alpha \in H^1(F, G)$ with zero Rost invariant, and write $\beta$ for its image in $H^1(F, \aut(\G(C, L)))$; it corresponds to a twisted composition.
The image of $\beta$ in $H^1(F, \aut(J(\G(C, L))))$ also has Rost invariant zero by \cite[p.~122]{GMS},\footnote{The inclusion $G \subset \aut(J(\G(C, L))$ has Rost multiplier 1, but we do not need this fact here.} hence is zero by the lemma.  In particular, the Albert $F$-algebra obtained by the Springer Decomposition from $\beta$ is has zero divisors and so $\beta$ is a twisted Hurwitz composition $\G(C', L')$ by \cite[38.8]{KMRT}.  But $\beta$ comes from $\alpha$, so in fact $L' = L$.  And reading the statement of \cite[38.8]{KMRT} closely, we see that $C'$ is the coordinate algebra of the split Albert algebra, i.e., $C'$ is the split octonion algebra $C$.  We have deduced that $\alpha$ is in the kernel of the map $H^1(F, G) \ra H^1(F, \aut(\G(C, L)))$.

Write $\Gb$ for the adjoint quotient of $G$ and $\D$ for its Dynkin diagram endowed with the natural action of the Galois group of $F$.  As $G$ is quasi-split, so is $\Gb$, and the sequence (from, say, \cite[16.3.9(3)]{Sp:LAG} or \cite[Lemma 1.6]{MPW1})
\[
1 \ra \Gb \ra \aut(\Gb) \ra \aut(\D) \ra 1
\]
is exact on $F$-points.  Said differently, the map $H^1(F, \Gb) \ra H^1(F, \aut(\Gb))$ has zero kernel.  But $\aut(\G(C, L))$ is just $G \rtimes \aut(\D)$ and $\aut(\Gb)$ is $\Gb \rtimes \aut(\D)$, so there is a commutative diagram
\[
\begin{CD}
H^1(F, G) @>>> H^1(F, \aut(\G(C, L))) \\
@VVV @VVV \\
H^1(F, \Gb) @>>> H^1(F, \aut(\Gb))
\end{CD}
\]
Clearly, $\alpha$ is in the kernel of $H^1(F, G) \ra H^1(F, \Gb)$.  This is the same as the image of $H^1(F, Z) \ra H^1(F, G)$, where $Z$ denotes the center of $G$.  As $G$ is quasi-split, it contains a quasi-trivial maximal torus $S$ and it in turn contains $Z$ (because every maximal torus does).  Then $\alpha$ is in the image of $H^1(F, S) \ra H^1(F, G)$, but $H^1(F, S)$ is zero by Hilbert's Theorem 90.
\end{proof}

\providecommand{\bysame}{\leavevmode\hbox to3em{\hrulefill}\thinspace}
\providecommand{\MR}{\relax\ifhmode\unskip\space\fi MR }
\providecommand{\MRhref}[2]{%
  \href{http://www.ams.org/mathscinet-getitem?mr=#1}{#2}
}
\providecommand{\href}[2]{#2}

\end{document}

%% file: skip-thm.tex
\usepackage{amsthm}

\newtheorem*{thm*}{Theorem}
\newtheorem*{prop*}{Proposition}
\newtheorem*{cor*}{Corollary}
\newtheorem*{lem*}{Lemma}
\newtheorem*{MT*}{Main Theorem}

\theoremstyle{definition} %

\newtheorem*{defn*}{Definition}

\theoremstyle{remark} %

\newtheorem*{rmk*}{Remark}
\newtheorem*{rmks*}{Remarks}

\newtheoremstyle{exercise}
  {3pt}
  {3pt}
  {\small}
  {}
  {\sc\small}
  {.}
  {.5em}
   {}     
  {}

\theoremstyle{exercise}

\makeatletter
{\renewcommand{\theequation}{#1}}%
{\renewcommand{\theequation}{\arabic{equation}}\addtocounter{equation}{-1}\global\@ignoretrue}
\makeatother

\makeatletter
{\renewcommand{\theequation}{#1}\begin{eqnarray}}%
{\end{eqnarray}\renewcommand{\theequation}{\arabic{equation}}\addtocounter{equation}{-1}\global\@ignoretrue}
\makeatother

\makeatletter
{\smallskip \refstepcounter{equation}\noindent{\textbf{\theequation.} }{{\textbf{#1.}}}}%
{\smallskip \global\@ignoretrue}
\makeatother

\makeatletter
{\smallskip \refstepcounter{equation}\noindent{\textbf{\theequation.} }{{\textbf{#1}}}}%
{\smallskip \global\@ignoretrue}
\makeatother

\makeatletter
{\smallskip \refstepcounter{equation}{\sc \theequation}{\sc (#1).}}%
{\smallskip \global\@ignoretrue}
\makeatother

\makeatletter
{\smallskip\refstepcounter{equation}\noindent{\textbf{\theequation.}}{\textsl{ #1.}}}%
{\smallskip\global\@ignoretrue}
\makeatother

\makeatletter
\newenvironment{borel*}%
{\smallskip \refstepcounter{equation}\noindent{\textbf{\theequation.}}}%
{\global\@ignoretrue}
\makeatother

\newcommand{\flist}[1]{\hangindent\leftmargini\textup{(1)}\hskip\labelsep {#1}%
\begin{enumerate}%
\setcounter{enumi}{1}%
}

%% file: qformsl-def.tex
\newcommand{\C}{{\mathbb{C}}}        
\newcommand{\Q}{{\mathbb{Q}}}        
\newcommand{\Z}{{\mathbb{Z}}}        

\newcommand{\QZt}{\mathbb{Q}/\mathbb{Z}(2)}

\newcommand{\Zm}[1]{\Z/{#1}\Z}

\newcommand{\D}{\Delta}

\newcommand{\G}{{\Gamma}}       

\newcommand{\oddots}{{\mathinner{\mkern1mu\raise1pt\vbox{\kern7pt\hbox{.}}\mkern2mu\raise4pt\hbox{.}\mkern2mu\raise7pt\hbox{.}\mkern1mu}}}

\newcommand{\iiiD}{^3\!D_4}
\newcommand{\viD}{^6\!D_4}

\DeclareMathOperator{\aut}{Aut}

\newcommand{\injects}{\hookrightarrow}

\newcommand{\ra}{\rightarrow}